\documentclass[12pt, leqno]{amsart}

\usepackage{texdraw,multicol,rotating}

\usepackage{amsmath}
\usepackage{amssymb}
\usepackage{enumerate}

\input{txdtools.tex}

\numberwithin{equation}{section}

\theoremstyle{plain}
\newtheorem{thm}{Theorem}[section]
\newtheorem{prop}[thm]{Proposition}
\newtheorem{lem}[thm]{Lemma}

\theoremstyle{definition}
\newtheorem{defi}[thm]{Definition}

\newtheorem{example}[thm]{Example}
\newtheorem{remark}[thm]{Remark}

\theoremstyle{remark}

\newbox{\tmpa}
\newbox{\tmpb}

\DeclareMathOperator{\wt}{wt}
\newcommand{\nc}{\newcommand}
\nc{\Uq}{U_q} \nc{\Z}{\mathbf{Z}} \nc{\C}{\mathbf{C}}
\nc{\Q}{\mathbf{Q}}
\nc{\op}{\oplus} \nc{\ot}{\otimes} \nc{\pv}{P^{\vee}}
\nc{\ali}{\alpha_i} \nc{\B}{\mathbf{B}} \nc{\F}{\mathbf{F}}
\nc{\bP}{\mathbf{P}} \nc{\V}{\mathbf{V}} \nc{\La}{\Lambda}
\nc{\la}{\lambda}
\nc{\nbinom}[2]{\genfrac{}{}{0pt}{1}{{#1}}{{#2}}}
\nc{\qbinom}[2]{\left[\genfrac{}{}{0pt}{1}{{#1}}{{#2}}\right]}
\nc{\path}{\mathcal{P}} \nc{\fit}{\tilde{f}_i}
\nc{\eit}{\tilde{e}_i} \nc{\fjt}{\tilde{f}_j}
\nc{\ejt}{\tilde{e}_j} \nc{\Y}{\mathbf{Y}} \nc{\A}{\mathbf{A}}
\nc{\ra}{\rightarrow} \nc{\vep}{\varepsilon} \nc{\vphi}{\varphi}
\nc{\vp}{\varphi} \nc{\g}{\mathfrak{g}} \nc{\h}{\mathfrak{h}}
\nc{\oP}{\overline{P}} \nc{\pathp}{\mathbf{p}}

\nc{\tris}{ \bsegment \move(0 0)\lvec(10 0)\lvec(10 10)\lvec(0
0)\ifill f:0.7 \esegment } \nc{\recs}{ \bsegment \move(0
0)\lvec(10 0)\lvec(10 5)\lvec(0 5)\lvec(0 0)\ifill f:0.7 \esegment
}
\nc{\hcvec}[5]{%
\getpos(#1 #3)\spx\spy \getpos(#2 #3)\epx\epy \getpos(#4
#5)\xoff\yoff \realadd \spx \xoff \twox \realadd \epx {-\xoff}
\thrx \realadd \spy \yoff \posy \move({\spx} {\spy}) \clvec
({\twox} {\posy})({\thrx} {\posy})({\epx} {\epy}) \rmove(0 0) }

\nc{\ahead}[2]{%
\cossin (0 0)({#1} {#2})\cosa\sina \bsegment
  \drawdim in \setunitscale 0.065
  \realmult {-0.5} \cosa \hcosa
  \realmult {-0.5} \sina \hsina
  \move({\hcosa} {\hsina}) \ravec({\cosa} {\sina})
\esegment }
\nc{\boxi}{%
{%
\savebox{\tmppic}{\begin{texdraw} \small \drawdim em \textref h:C
v:C \setunitscale 0.55 \htext(0 0){$i$} \move(-1 -1)\lvec(-1
1)\lvec(1 1)\lvec(1 -1)\lvec(-1 -1)
\end{texdraw}}%
\raisebox{-0.19\height}{\usebox{\tmppic}}%
}%
}
\nc{\boxj}{%
{%
\savebox{\tmppic}{\begin{texdraw} \small \drawdim em \textref h:C
v:C \setunitscale 0.55 \htext(0 0.1){$j$} \move(-1 -1)\lvec(-1
1)\lvec(1 1)\lvec(1 -1)\lvec(-1 -1)
\end{texdraw}}%
\raisebox{-0.19\height}{\usebox{\tmppic}}%
}%
}
\nc{\boxipo}{%
{%
\savebox{\tmppic}{\begin{texdraw} \small \drawdim em \textref h:C
v:C \setunitscale 0.55 \htext(0.15 0){$i\!\!+\!\!1$} \move(-1.4
-1)\lvec(-1.4 1)\lvec(1.4 1)\lvec(1.4 -1)\lvec(-1.4 -1)
\end{texdraw}}%
\raisebox{-0.19\height}{\usebox{\tmppic}}%
}%
} \everytexdraw{ \drawdim in \arrowheadsize l:0.065 w:0.03
\arrowheadtype t:F \textref h:C v:C }
\newsavebox{\tmppic}
\newsavebox{\tmpfig}
\newsavebox{\tmpdraw}
\newsavebox{\tmpfiga}
\newsavebox{\tmpfigb}
\newsavebox{\tmpfigc}
\newsavebox{\tmpfigd}
\newsavebox{\tmpfige}
\newsavebox{\tmpfigf}
\newsavebox{\tmpfigg}
\newsavebox{\tmpfigh}
\newsavebox{\tmpfigi}
\newsavebox{\tmpfigj}
\newsavebox{\tmpfigk}
\newsavebox{\tmpfigl}
\newsavebox{\tmpfigm}
\newsavebox{\tmpfign}
\newsavebox{\tmpfigo}
\newsavebox{\tmpfigp}
\newsavebox{\tmpfigq}
\newsavebox{\tmpfigr}
\newsavebox{\tmpfigs}
\newsavebox{\tmpfigt}
\newsavebox{\tmpfigu}
\newsavebox{\tmpfigv}
\newsavebox{\tmpfigw}
\newsavebox{\tmpfigx}
\newsavebox{\tmpfigy}
\newsavebox{\tmpfigz}

\newsavebox{\tmpfigaa}
\newsavebox{\tmpfigab}
\newsavebox{\tmpfigac}
\newsavebox{\tmpfigad}
\newsavebox{\tmpfigae}
\newsavebox{\tmpfigaf}
\newsavebox{\tmpfigag}
\newsavebox{\tmpfigah}
\newsavebox{\tmpfigai}
\newsavebox{\tmpfigaj}
\newsavebox{\tmpfigak}
\newsavebox{\tmpfigal}
\newsavebox{\tmpfigam}
\newsavebox{\tmpfigan}
\newsavebox{\tmpfigao}
\newsavebox{\tmpfigap}
\newsavebox{\tmpfigaq}
\newsavebox{\tmpfigar}
\newsavebox{\tmpfigas}
\newsavebox{\tmpfigat}
\newsavebox{\tmpfigau}
\newsavebox{\tmpfigav}
\newsavebox{\tmpfigaw}
\newsavebox{\tmpfigax}
\newsavebox{\tmpfigay}
\newsavebox{\tmpfigaz}

\newsavebox{\tmpfigba}
\newsavebox{\tmpfigbb}
\newsavebox{\tmpfigbc}
\newsavebox{\tmpfigbd}
\newsavebox{\tmpfigbe}
\newsavebox{\tmpfigbf}
\newsavebox{\tmpfigbg}
\newsavebox{\tmpfigbh}

\nc{\node}{\lcir r:1 }
\nc{\sline}{\bsegment\savepos(10 0)(*ex *ey)
            \move(1 0)\rlvec(8 0)
            \esegment\move(*ex *ey)}
\nc{\dline}{\bsegment\savepos(10 0)(*ex *ey)
            \move(0.93 0.4)\rlvec(8.14 0)\rmove(0 -0.8)\rlvec(-8.14 0)
            \esegment\move(*ex *ey)}
\nc{\uline}{\bsegment\savepos(0 10)(*ex *ey)
            \move(0 1)\rlvec(0 8)
            \esegment\move(*ex *ey)}
\nc{\lpoint}{\savecurrpos(*ex *ey)
             \rmove(2.5 1.5)\rlvec(-1.5 -1.5)\rlvec(1.5 -1.5)
             \move(*ex *ey)}
\nc{\rpoint}{\savecurrpos(*ex *ey)
             \rmove(-2.5 -1.5)\rlvec(1.5 1.5)\rlvec(-1.5 1.5)
             \move(*ex *ey)}
\nc{\bline}{\bsegment\savepos(10 0)(*ex *ey)
            \linewd 0.6 \move(1.1 0)\rlvec(7.8 0)
            \esegment\move(*ex *ey)}

\nc{\araise}[1]{\raisebox{4.5pt}{#1}}
\nc{\braise}[1]{\raisebox{12.1pt}{#1}}
\nc{\craise}[1]{\raisebox{8pt}{#1}}
\nc{\draise}[1]{\raisebox{12pt}{#1}} \nc{\be}{\begin{enumerate}}
\nc{\ee}{\end{enumerate}} \nc{\bnum}{\begin{enumerate}[{\rm(i)}]}
\nc{\cl}{\colon} \nc{\seteq}{\mathbin{:=}} \nc{\re}{\mathrm{re}}
\nc{\im}{\mathrm{im}} \nc{\ran}{\rangle} \nc{\lan}{\langle}
\nc{\on}{\operatorname}
\newcommand{\set}[2]{\left\{#1\,;\,#2\,\right\}}
\nc{\Hom}{\on{Hom}} \nc{\Oint}{\mathcal{O}_{\mathrm{int}}}
\nc{\Wt}{\on{Wt}} \nc{\pP}{\widetilde{P}}
\nc{\eq}{\begin{eqnarray}} \nc{\eneq}{\end{eqnarray}}
\nc{\eqn}{\begin{eqnarray*}} \nc{\eneqn}{\end{eqnarray*}}
\nc{\Lemma}{\begin{lem}} \nc{\enlemma}{\end{lem}}
\newcommand{\isoto}[1][]{\mathop{\xrightarrow[#1]%
{\rule{0pt}{.9ex}%
{\raisebox{-.35ex}[0ex][-.6ex]{$\mspace{1mu}\sim\mspace{2mu}$}}}}}

\newcommand{\To}[1][\phantom{aaaa}]{\xrightarrow{\;#1\;}}
\nc{\hs}{\hspace*} \nc{\bfi}{\mathbf{i}} \nc{\eps}{\varepsilon}
\nc{\ba}{\begin{array}} \nc{\ea}{\end{array}}
\renewcommand{\Im}{\operatorname{Im}}
\nc{\tf}{\tilde{f}} \nc{\id}{\operatorname{id}} \nc{\bl}{\bigl}
\nc{\br}{\bigr} 

\begin{document}

\title[Abstract Crystals for Generalized Kac-Moody Algebras]
      {Abstract Crystals for Quantum \\ Generalized Kac-Moody Algebras}
\author[K. Jeong, S.-J. Kang, M. Kashiwara,  D.-U. Shin]
{Kyeonghoon Jeong$^{*}$, Seok-Jin Kang$^{*}$, Masaki
Kashiwara$^{\dagger}$ and\\ Dong-Uy Shin$^{\diamond}$}
\address{$^{*}$ Department of Mathematical Sciences and Research
Institute of Mathematics\\Seoul National University \\
San 56-1 Shinrim-Dong, Kwanak-Ku\\ Seoul, 151-747, Korea}
\email{khjeong@math.snu.ac.kr, sjkang@math.snu.ac.kr}

\address{$^{\dagger}$ Research Institute for Mathematical Sciences, Kyoto
University, Kitashirakawa, Sakyo-Ku, Kyoto 606--8502, Japan}
\email{masaki@kurims.kyoto-u.ac.jp}

\address{$^{\diamond}$ Department of Mathematics, Chonnam National
University, Yongbong-dong, Buk-gu\\ Kwangju, 500-757, Korea}
\email{dushin@chonnam.ac.kr}

\thanks{$^{*}$ This research was supported in part by KOSEF Grant
\# R01-2003-000-10012-0 and KRF Grant \# 2005-070-C00004. \hfill \\
$^{\diamond}$ This research was supported by KOSEF Grant
\# R01-2003-000-10012-0. \hfill \\
AMS Classification 2000: 17B37, 81R50}

\begin{abstract}

In this paper, we introduce the notion of abstract crystals for
quantum generalized Kac-Moody algebras and study their fundamental
properties. We then prove the crystal embedding theorem and give a
characterization of the crystals $B(\infty)$ and $B(\la)$.

\end{abstract}

\maketitle

\vskip 1cm \setcounter{page}{1}
\section*{Introduction}

The purpose of this paper is to develop the theory of {\it abstract
crystals} for quantum generalized Kac-Moody algebras. In
\cite{Kas91}, the third author introduced the {\it crystal basis
theory} for quantum groups associated with symmetrizable Kac-Moody
algebras. (In \cite{Lusztig90}, Lusztig constructed canonical bases
for quantum groups of ADE type.)  It has become one of the most
central themes in combinatorial representation theory, for it
provides us with a very powerful combinatorial tool to investigate
the structure of integrable modules over quantum groups and
Kac-Moody algebras.

\vskip 2mm

The {\it generalized Kac-Moody algebras} were introduced by
Borcherds in his study of Monstrous Moonshine \cite{Bo88}. The
{\it Monster Lie algebra}, an example of generalized Kac-Moody
algebras, played a crucial role in his proof of the Moonshine
conjecture \cite{Bo92}. In \cite{Kang95}, the second author
constructed the {\it quantum generalized Kac-Moody algebra}
$U_q(\g)$ as a deformation of the universal enveloping algebra of
a generalized Kac-Moody algebra $\g$. He also showed that, for a
generic $q$, the Verma modules and the unitarizable highest weight
modules over $\g$ can be deformed to those over $U_q(\g)$ in such
a way that the dimensions of weight spaces are invariant under the
deformation.

\vskip 2mm

In \cite{JKK}, the first three authors developed the crystal basis
theory for quantum generalized Kac-Moody algebras. More precisely,
they defined the notion of crystal bases for $U_q(\g)$-modules in
the category $\Oint$ (see \S\;\ref{sec:gen}), proved standard
properties of crystal bases including the tensor product rule, and
showed that there exists a crystal basis (and a global basis) of
the negative part $U_q^{-}(\g)$ of a quantum generalized Kac-Moody
algebra and one of the irreducible $U_q(\g)$-module $V(\la)$ with
a dominant integral weight $\la$ as its highest weight.

\vskip 2mm

In this paper, we introduce the notion of abstract crystals for
quantum generalized Kac-Moody algebras and investigate their
fundamental properties.
We then prove the {\it crystal embedding theorem}, which yields a
procedure to determine the structure of the crystal $B(\infty)$ in
terms of elementary crystals. Finally, as an application of the
crystal embedding theorem, we provide a characterization of the
crystals $B(\infty)$ and $B(\la)$.  We also include an explicit
description of the crystals $B(\infty)$ and $B(\la)$ for quantum
generalized Kac-Moody algebras of rank 2 and for the quantum
Monster algebra.

\section{Generalized Kac-Moody algebras}
\label{sec:gen}

Let $I$ be a finite or countably infinite index set. A real matrix
$A=(a_{ij})_{i, j\in I}$ is called a {\it Borcherds-Cartan matrix}
if it satisfies the following conditions: \bnum   \item $a_{ii}=2$
or $a_{ii} \le 0$ for all $i\in I$, \item $a_{ij} \le 0$ if $i
\neq j$, \item $a_{ij} \in \Z$ if $a_{ii}=2$, \item $a_{ij}=0$ if
and only if $a_{ji}=0$. \ee In this paper, {\em we assume that $A$
is even and integral}\,; i.e., $a_{ii} \in 2\Z_{\le1}$ for all
$i\in I$ and $a_{ij} \in \Z$ for all $i, j \in I$. Furthermore, we
also assume that $A$ is {\it symmetrizable}; i.e., there exists a
diagonal matrix $D=\text{diag} (s_i \in \Z_{>0} ; i \in I)$ such
that $DA$ is symmetric.

We say that an index $i\in I$ is {\it real} if $a_{ii}=2$ and {\it
imaginary} if $a_{ii}\le 0$. We denote by $I^\re=\set{i\in
I}{a_{ii}=2}$ and $I^\im=\set{ i \in I}{a_{ii} \le 0 }$ the set of
real indices and the set of imaginary indices, respectively.

A {\it Borcherds-Cartan datum} $(A, P, \Pi, \Pi^{\vee})$ consists
of

\bnum
\item a Borcherds-Cartan matrix $A=(a_{ij})_{i, j \in I}$,

\item a free abelian group $P$, the {\em weight lattice},

\item $\Pi=\set{\alpha_i\in P}{i\in I}$, the set of {\em
simple roots},

\item $\Pi^{\vee} = \set{ h_i}{ i\in I }\subset
P^\vee\seteq\Hom(P,\Z)$, the set of {\em simple coroots}, \ee
satisfying the properties: \be[\quad(a)] \item $\langle
h_i,\alpha_j \rangle = a_{ij}$ for all $i, j \in I$,

\item for any $i\in I$, there exists $\Lambda_i\in P$ such that
$\lan h_j,\Lambda_i\ran=\delta_{ij}$ for all $j\in I$,

\item $\Pi$ is linearly independent.
\ee
We denote by $P^{+}=\set{\la \in P}{ \la(h_i) \ge 0 \ \ \text{for
all} \ i\in I }$ the set of {\it dominant integral weights}. We
also use the notation $Q=\bigoplus_{i\in I} \Z \alpha_i$ and
$Q_{+}=\sum_{i\in I} \Z_{\ge 0} \alpha_i$.

Let $q$ be an indeterminate and set $q_i = q^{s_i}$ $(i\in I)$.
For an integer $n \in \Z$, define
\begin{equation*}
[n]_i = \dfrac{q_i^n - q_i^{-n}} {q_i - q_i^{-1}}, \qquad [n]_i! =
\prod_{k=1}^n [k]_i, \qquad \left[\begin{matrix} m \\
n \end{matrix}\right]_i 
= \dfrac{[m]_i!}{[n]_i! [m-n]_i!}.
\end{equation*}

Let $(A, P, \Pi, \Pi^{\vee})$ be a Borcherds-Cartan datum. The
{\it quantum generalized Kac-Moody algebra} $U_q(\g)$ associated
with $(A, P, \Pi,\Pi^{\vee})$ is defined to be the associated
algebra over $\Q(q)$ with 1 generated by the elements $e_i$, $f_i$
$(i\in I)$, $q^h$ $(h\in P^{\vee})$ with the following defining
relations: {\allowdisplaybreaks
\begin{equation} 
\begin{aligned}
& q^0 =1, \quad q^h q^{h'} = q^{h+h'}\quad\text{for} \ \ h, h' \in
P^{\vee}, \\
& q^h e_i q^{-h} = q^{\alpha_i(h)} e_i, \quad q^h f_i q^{-h} =
q^{-\alpha_i(h)} f_i \quad\text{for $h\in P^{\vee}$, $i\in I$,} \\
& e_i f_j - f_j e_i = \delta_{ij} \dfrac{K_i - K_i^{-1}} {q_i -
q_i^{-1}} \quad\text{for $i, j \in I$, where $ K_i =
q^{s_i h_i}$,} \\
& \sum_{k=0}^{1-a_{ij}} (-1)^k \left[\begin{matrix} 1-a_{ij} \\ k
\end{matrix} \right]_{i} e_i^{1-a_{ij}-k} e_j e_i^k =0 \quad\text{if
$i\in I^\re$ and $i\neq j$,}\\
& \sum_{k=0}^{1-a_{ij}} (-1)^k \left[\begin{matrix} 1-a_{ij} \\ k
\end{matrix} \right]_{i} f_i^{1-a_{ij}-k} f_j f_i^k =0 \quad
\text{if $i\in I^\re$ and $i\neq j$,} \\
& e_i e_j - e_j e_i = f_i f_j - f_j f_i =0 \quad\text{if
$a_{ij}=0$.}
\end{aligned}
\end{equation}
} Let us denote by $U_q^+(\g)$ (resp.\ $U_q^-(\g)$) the subalgebra
of $U_q(\g)$ generated by the $e_i$'s (resp.\ the $f_i$'s). Let us
denote by $\Oint$ the abelian category of $U_q(\g)$-modules $M$
satisfying the following properties: \bnum \item $M$ has the
weight decomposition: $M=\oplus_{\lambda\in P}M_\lambda$, where
$$M_\la\seteq \set{u\in M}{q^hu=q^{\la(h)}u\quad\text{for any
$h\in P^\vee$}},$$ \item the action of $U_q^+(\g)$ is locally
finite, i.e., $\dim U_q^+(\g)u<\infty$ for any $u\in M$, \item
$\wt(M)\seteq\set{\la\in P}{M_\la\not=0} \subset \set{\la\in
P}{\text{$\lan h_i,\la\ran\ge0$ for any $i\in I^\im$}}$, \item
$f_iM_\la=0$  for any $i\in I^\im$ and $\la\in P$ such that $\lan
h_i,\la\ran=0$, \item $e_iM_\la=0$  for any $i\in I^\im$ and
$\la\in P$ such that $\lan h_i,\la\ran\le-a_{ii}$. \ee It is
proved in \cite{JKK} that the abelian category $\Oint$ is
semisimple, and any of its irreducible objects is of the form
$V(\la)$ for $\la\in P^+$, where $V(\la)=U_q(\g)u_\la$ with the
defining relation: \be[{\quad(a)}] \item $u_\la$ has weight $\la$,
\item $e_iu_\la=0$ for all $i\in I$,

\item $f_i^{\lan h_i,\la\ran+1}u_\la=0$ for any $i\in I^\re$,

\item
$f_iu_\la=0$ if $i\in I^\im$ and $\lan h_i,\la\ran=0$. \ee

Let $\A_0 = \set{f/g \in \Q(q) }{ f, g \in \Q[q], g(0) \neq 0 }$.
Let $M$ be a $U_q(\g)$-module in $\Oint$. For each $i\in I$, every
weight vector $u\in M_{\mu}$ has an {\it $i$-string decomposition}
\begin{equation*}
u= \sum_{k\ge 0} f_i^{(k)} u_k\quad\text{with $u_k\in
M_{\mu+n\alpha_i} $ such that $e_i u_k =0$,}
\end{equation*}
where
\begin{equation*}
f_i^{(k)} = \begin{cases} f_i^k / [k]_i! \quad & \text{if} \ \
i\in I^\re, \\
f_i^k \quad & \text{if} \ \ i\in I^\im.
\end{cases}
\end{equation*}
Such a decomposition is unique. We define the {\it Kashiwara
operators} $\eit$, $\fit$ $(i\in I)$ by
\begin{equation*}
\eit u = \sum_{k\ge 1} f_i^{(k-1)} u_k, \qquad \fit u = \sum_{k\ge
0} f_i^{(k+1)} u_k.
\end{equation*}

A crystal basis $(L,B)$ of $M$ is a pair of a free
$\A_0$-submodule $L$ of $M$ and a basis $B$ of the $\Q$-vector
space $L/qL$ satisfying the following conditions, \bnum
\item $L$ generates $M$ as a $\Q(q)$-vector space,
\item $L$ has the weight decomposition $L=\oplus_{\la\in P}L_\la$
where $L_\la=L\cap M_\la$,
\item $B$ has the weight decomposition $B=\sqcup_{\la\in P}B_\la$
where $B_\la=B\cap (L_\la/qL_\la)$,
\item $\fit L\subset L$ and $\eit L\subset L$ for any $i\in I$,
\item $\fit B\subset B\cup\{0\}$ and $\eit B\subset B\cup\{0\}$
for any $i\in I$,
\item for $b,b'\in B$ and $i\in I$, $\fit b =b'$ if and only if $b=\eit b'$,
\ee It is proved in \cite{JKK} that every $M\in\Oint$ has a
crystal basis unique up to an automorphism.

For $\la\in P^+$, let $L(\la)$ be the $\A_0$-submodule of
$V(\lambda)$ generated by
$$\set{\tilde f_{i_1} \cdots \tilde f_{i_r}
u_{\la}}{r \ge 0, i_k \in I },$$ and set $B(\la) = \set{\tilde
f_{i_1} \cdots \tilde f_{i_r} u_{\la} + q L(\la) }{ r\ge 0, i_k
\in I } \setminus \{0\}$. Then $(L(\lambda),B(\la))$ is a crystal
basis of $V(\lambda)$.

\section{Abstract Crystals}
By abstracting the properties of crystal bases of
$U_q(\g)$-modules in $\Oint$, we shall introduce the notion of
abstract crystals.

\begin{defi} \label{defi:abstract crystal}
An {\em abstract $U_q(\g)$-crystal} or simply a {\em crystal\/} is
a set $B$ together with the maps $\wt\cl B \ra P$, $\eit, \fit\cl
B \ra B \sqcup \{0\}$ and $\vep_i, \vphi_i\cl  B \ra \Z \sqcup
\{-\infty\}$ $(i\in I)$ satisfying the following conditions: \bnum
\item $\wt(\eit b) = \wt b + \alpha_i$ if $i\in I$ and
$\eit b \neq 0$,

\item $\wt(\fit b) = \wt b - \alpha_i$ if $i\in I$ and
$\fit b \neq 0$,

\item for any $i \in I$ and $b\in B$, $\vphi_i(b) = \vep_i(b) +
\langle h_i, \wt b \rangle$,

\item for any $i\in I$ and $b,b'\in B$,
$\fit b = b'$ if and only if $b = \eit b'$,\label{cond7}

\item for any $i \in I$ and $b\in B$
such that $\eit b \neq 0$, we have \be[{\rm(a)}]
\item
$\vep_i(\eit b) = \vep_i(b) - 1$, $\vphi_i(\eit b) = \vphi_i(b) +
1$ if $i\in I^\re$,
\item
$\vep_i(\eit b) = \vep_i(b)$ and $\vphi_i(\eit b) = \vphi_i(b) +
a_{ii}$ if $i\in I^\im$, \ee

\item for any $i \in I$ and $b\in B$ such that $\fit b \neq 0$,
we have \be[{\rm(a)}]
\item
$\vep_i(\fit b) = \vep_i(b) + 1$ and $\vphi_i(\fit b) = \vphi_i(b)
- 1$ if $i\in I^\re$,
\item
$\vep_i(\fit b) = \vep_i(b)$ and $\vphi_i(\fit b) = \vphi_i(b) -
a_{ii}$ if $i\in I^\im$, \ee
\item for any $i \in I$ and $b\in B$ such that $\vphi_i(b) = -\infty$, we
have $\eit b = \fit b = 0$.

\ee
\end{defi}
We sometimes write
$$\wt_i(b)=\lan h_i,\wt b\ran\quad\text{for $i\in I$ and $b\in B$.}$$

\begin{remark}
Almost all crystals appearing in this paper have the following
properties: \be[\quad(a)]
\item $\wt_i(b)\ge0$ for any $i\in I^\im$ and $b\in B$,
\item $\vep_i(b)\in \Z_{\le0}\sqcup\{-\infty\}$
and $\vp_i(b)\in \Z_{\ge0}\sqcup\{-\infty\}$ for any $i\in I^\im$
and $b\in B$, \ee Hence we could include these properties in the
axiom of crystals.
\end{remark}

We shall define the morphisms of crystals.
\begin{defi}\label{def:mor}
Let $B_1$ and $B_2$ be crystals. A {\em morphism of crystals} or a
{\em crystal morphism} $\psi\cl B_1\rightarrow B_2$ is a map
$\psi\cl B_1\to B_2$ such that

\bnum
\item for $b\in B_1$ we have
\begin{equation*}
\text{$\wt(\psi(b))=\wt(b)$, and $\vep_i(\psi(b))=\vep_i(b)$,
$\vp_i(\psi(b))=\vp_i(b)$ for all $i\in I$,}
\end{equation*}

\item if $b\in B_1$ and $i\in I$ satisfy $\fit b\in B_1$, then we have
$\psi(\fit b)=\fit\psi(b)$. \label{cond:crysmor2} \ee
\end{defi}

\begin{remark}
If $b\in B_1$ satisfies $\eit b\in B_1$, one can deduce $\psi(\eit
b)=\eit\psi(b)$ for all $i\in I$ using Definition~
\ref{defi:abstract crystal} \eqref{cond7} and
Definition~\ref{def:mor} \eqref{cond:crysmor2}.
\end{remark}

Then the crystals form a category.

\begin{defi}
Let $\psi\cl B_1\rightarrow B_2$ be a morphism of crystals.
\be[{\quad \rm(a)}]
\item $\psi$ is called a {\it strict morphism} if
\begin{equation*}
\psi(\eit b)=\eit\psi(b),\,\, \psi( \fit
b)=\fit\psi(b)\quad\text{for all $i\in I$ and $b\in B_1$.}
\end{equation*}
Here we understand $\psi(0)=0$.
\item $\psi$ is called an {\em embedding} if the underlying
map $\psi\cl B_1\rightarrow B_2$ is injective. In this case, we
say that $B_1$ is a {\em subcrystal} of $B_2$. If $\psi$ is a
strict embedding, we say that $B_1$ is a {\em full subcrystal} of
$B_2$.


\ee
\end{defi}

\begin{remark}

If $B_1$ is a full subcrystal of $B_2$, then we have
\begin{equation*}
B_2\cong B_1\oplus(B_2\setminus B_1).
\end{equation*}
\end{remark}

\medskip
Let us give two examples of crystals.

\begin{example}[\cite{JKK}]
Let $(L,B)$ be a crystal basis of $M\in\Oint$. Then, $B$ has a
crystal structure, where the maps $\vep_i, \vp_i$ ($i\in I$) are
given by
\begin{equation*}
\begin{aligned}
\vep_i (b) & =
\begin{cases}
\max \set{k\ge 0 }{ \eit^k b \neq 0 }&\text{for $i\in I^\re$,} \\
0&\text{for $i\in I^\im$,}\end{cases}\\
\vphi_i (b) & =
\begin{cases}\max \set{k\ge 0 }{ \fit^k b \neq 0 }&\text{for $i\in I^\re$,}
\\
\wt_i(b)&\text{for $i\in I^\im$,}\end{cases}
\end{aligned}
\end{equation*}
(This definition is different from the one given in \cite{JKK}.)
Such a crystal $B$ has the following properties: \be[\quad(a)]
\item $\vep_i(b),\vp_i(b)\ge0$
for any $b\in B$ and $i\in I$,
\item if $i\in I^\im$ and $\wt_i(b)=0$, then $\eit b=\fit b=0$,
\item if $i\in I$ and $\vp_i(b)>0$, then $\fit b\not=0$.
\ee

Let $V(\la)$ be the irreducible highest weight $U_q(\g)$-module
with highest weight $\la \in P^{+}$. Let us recall that $V(\la)$
has a crystal basis $(L(\lambda),B(\la))$. Hence $B(\la)$ has a
crystal structure.

\end{example}

\begin{example}[\cite{JKK}]
Fix $i\in I$. For any $u \in U_q^{-}(\g)$, there exist unique $v,
w \in U_q^{-}(\g)$ such that
\begin{equation*}
e_i u - u e_i = \dfrac{K_i v - K_i^{-1} w}{q_i - q_i^{-1}}.
\end{equation*}
We define the endomorphism $ e_i' \cl  U_q^{-}(\g) \ra
U_q^{-}(\g)$ by $e_i'(u)=w$. Then every $u\in U_q^{-}(\g)$ has a
unique $i$-string decomposition
\begin{equation*}
u=\sum_{k\ge 0} f_i^{(k)} u_k, \quad \text{where} \ \ e_i' u_k =0
\ \ \text{for all} \ \ k\ge 0,
\end{equation*}
and the Kashiwara operators $\eit$, $\fit$ $(i\in I)$ are defined
by
\begin{equation*}
\eit u = \sum_{k\ge 1} f_i^{(k-1)} u_k, \qquad \fit u = \sum_{k\ge
0} f_i^{(k+1)} u_k.
\end{equation*}
Let $L(\infty)$ be the $\A_0$-submodule of $U_q^{-}(\g)$ generated
by $$\set{\tilde f_{i_1} \cdots \tilde f_{i_r} \mathbf{1} }{ r \ge
0, i_k \in I },$$  and
\begin{equation*}
B(\infty) = \set{\tilde f_{i_1} \cdots \tilde f_{i_r} \mathbf{1} +
q L(\infty) }{ r\ge 0, i_k \in I } \setminus \{0\}\subset
L(\infty)/qL(\infty),
\end{equation*}
where $\mathbf{1}$ is the multiplicative identity in
$U_q^{-}(\g)$. Then $B(\infty)$ becomes a crystal with the maps
$\wt$, $\eit, \fit$, $\vep_i, \vp_i$ ($i\in I$), where
\begin{equation*}
\begin{aligned}
\wt(b) & =- (\alpha_{i_1} + \cdots + \alpha_{i_r}) \quad
\text{for} \ \ b=\tilde f_{i_1} \cdots \tilde f_{i_r}
\mathbf{1} + q L(\infty), \\
\vep_i (b) & =
\begin{cases} \max \set{k\ge 0 }{ \eit^k b \neq 0 }&\text{for $i\in I^\re$,}
\\
0&\text{for $i\in I^\im$,}\end{cases}\\
\vphi_i (b) & = \vep_i(b) + \wt_i(b) \quad (i\in I).
\end{aligned}
\end{equation*}
We have $\fit b\in B(\infty)$ for any $i\in I$ and $b\in
B(\infty)$.
\end{example}

The crystals $B(\la)$ and $B(\infty)$ are closely related as seen
in the following proposition.

\begin{prop}[\cite{JKK}] \label{prop:projection}
For every $\la \in P^{+}$, there exists a map $\pi_{\la}\cl
B(\la)\to B(\infty)$ such that \bnum
\item $\pi_\la$ is injective,
\item $\pi_{\la}(u_{\la}) = \mathbf{1}$,

\item $\pi_{\la} \circ \fit (b) = \fit \circ \pi_{\la} (b)$
for any $i\in I$ and $b\in B(\la)$ such that $\fit b\not=0$,

\item $\pi_{\la} \circ \eit (b) = \eit \circ \pi_{\la} (b)$
for all $i\in I$ and $b\in B(\la)$.

\item $\wt(\pi_{\la}(b))=\wt(b)-\la$,
$\vep_i(\pi_{\la}(b))=\vep_i(b)$ for any $b\in B(\la)$ and $i\in
I$. \ee
\end{prop}

\begin{proof}
See Propositions 7.12, 7.13, 7.23, 7.34 in \cite{JKK}.
\end{proof}

\bigskip
We define the tensor product of a pair of crystals as follows: for
two crystals $B_1$ and $B_2$, their tensor product $B_1\otimes
B_2$ is $\set{b_1\otimes b_2}{b_1\in B_1, b_2\in B_2}$ with the
following crystal structure. The maps $\wt, \vep_i,\vp_i$ are
given by \eqn
\wt(b\otimes b')&=&\wt(b)+\wt(b'),\\
\vep_i(b\otimes b')&=&\max(\vep_i(b),
\vep_i(b')-\wt_i(b)),\\
\vp_i(b\otimes b')&=&\max(\vp_i(b)+\wt_i(b'),\vp_i(b')). \eneqn
For $i\in I$, we define \eqn \fit(b\otimes b')&=&
\begin{cases}\fit b\otimes b'
&\text{if $\vp_i(b)>\vep_i(b')$,}\\
b\otimes \fit b' &\text{if $\vp_i(b)\le \vep_i(b')$,}
\end{cases}
\eneqn For $i\in I^\re$, we define \eqn \eit(b\otimes b')&=&
\begin{cases}\eit b\otimes b'\ &\text{if
$\vp_i(b)\ge \vep_i(b')$,}\\
b\otimes \eit b' &\text{if $\vp_i(b)< \vep_i(b')$,}
\end{cases}
\eneqn and, for $i\in I^\im$, we define \eqn \eit(b\otimes b')&=&
\begin{cases}\eit b\otimes b'\
&\text{if $\vp_i(b)>\vep_i(b')-a_{ii}$,}\\
0&\text{if $\vep_i(b')<\vp_i(b)\le\vep_i(b') -a_{ii}$,}\\
b\otimes \eit b' &\text{if $\vp_i(b)\le\vep_i(b')$.}
\end{cases}
\eneqn

(This tensor product rule is different from the one given in
\cite{JKK}. But when $B_1=B(\la)$ and $B_2 = B(\mu)$ for $\la, \mu
\in P^{+}$, the two rules coincide.)

\Lemma By the definition above,  $B_1\otimes B_2$ is a crystal.
\enlemma
\begin{proof}
The properties (i), (ii), (iii), (vii) are clear.

Suppose $\fit (b \ot b') = b_1 \ot b_1' \ne 0$ for $i \in
I^{\im}$. If $\vp_i(b) > \vep_i(b')$, then $b_1 = \fit b \ne 0$
and $b_1' = b'$. In this case, we have $\vp_i(b_1) = \vp_i (\fit
b) = \vp_i (b) - a_{ii}
>\vep_i(b')- a_{ii} = \vep_i(b_1') - a_{ii}$.
Hence we obtain $\eit (b_1 \ot b_1') = \eit b_1 \ot b_1' = b \ot
b'$. Also we have $\vep_i(b_1)=\vep_i(b) > \vep_i(b') - \wt_i(b)
\ge \vep_i(b') - \wt_i(\fit b)$. Hence $\vep_i (\fit (b \ot b')) =
\vep_i(b_1) = \vep_i (b \ot b')$.

If $\vp_i(b) \le \vep_i (b')$, then $b_1 = b$ and $b_1' = \fit b'
\ne 0$. In this case, we have $\vp_i(b_1) = \vp_i (b) \le
\vep_i(b') = \vep_i(b_1')$. Hence we obtain $\eit (b_1 \ot b_1') =
b_1 \ot \eit b_1' = b \ot b'$. Also we have $\vep_i(b) \le
\vep_i(b') - \wt_i(b) = \vep_i(\fit b') - \wt_i(b)$. Hence $\vep_i
(\fit (b \ot b')) = \vep_i(b') - \wt_i(b) = \vep_i (b \ot b')$.

Therefore (vi)(b) and the half of (iv) are proved. Now (v)(b)
follows easily from the property (iii). The rest may be proved
similarly. 
\end{proof}

Remark that, for a crystal basis $(L_i,B_i)$ of $M_i\in\Oint$
($i=1,2$),
$(L_1,B_1)\otimes(L_2,B_2)\seteq(L_1\otimes_{\A_0}L_2,B_1\otimes
B_2)$ is a crystal basis of $M_1\otimes M_2$, and the crystal
structure on $B_1\otimes B_2$ coincides with the tensor product of
the crystals $B_1$ and $B_2$.

It is also easy to check the following associativity law for the
tensor product, and the category of crystals has a structure of a
tensor category (see e.g.\ \cite{KS}). \Lemma For three crystals,
$B_\nu$ $(\nu=1,2,3)$, the map $(b_1\otimes b_2)\otimes b_3\mapsto
b_1\otimes (b_2\otimes b_3)$ gives an isomorphism of crystals:
$$\Phi : (B_1\otimes B_2)\otimes B_3\isoto B_1\otimes(B_2\otimes B_3).$$
\enlemma
\begin{proof}
We shall only show here that $\Phi(\eit ((b_1 \ot b_2) \ot b_3)) =
\eit (b_1 \ot (b_2 \ot b_3))$ for $i \in I^{\im}$, leaving the
proof of the rest to the reader. 

{\bf Case 1:} $\vp_i (b_1 \ot b_2) > \vep_i(b_3) - a_{ii}$.

In this case, we have $\eit ((b_1 \ot b_2) \ot b_3)=\bl(\eit(b_1
\ot b_2)\br)\otimes b_3$. 

If $\vp_i (b_1) > \vep_i (b_2) - a_{ii}$, then $\eit ((b_1 \ot
b_2) \ot b_3)=(\eit b_1 \ot b_2)\otimes b_3$. Since $\vp_i (b_1) >
\vep_i (b_2)$, we have $\vp_i(b_1 \ot b_2) = \vp_i (b_1) +
\wt_i(b_2)$. Hence we obtain $\vp_i(b_1) > \vep_i (b_3) - \wt_i
(b_2) - a_{ii}$. Since $\vp_i (b_1) > \vep_i (b_2) - a_{ii}$, we
obtain $\vp_i (b_1) > \vep_i (b_2 \ot b_3) - a_{ii}$, which
implies $\eit\bl(b_1 \ot( b_2 \ot b_3)\br)= \eit b_1\otimes( b_2
\ot b_3)$.


If $\vep_i (b_2) < \vp_i (b_1) \le \vep_i (b_2) - a_{ii}$, we will
show $\eit(b_1 \ot (b_2 \ot b_3)) = 0$. Since $\vep_i (b_2) <
\vp_i (b_1)$, we have $\vp_i (b_1) + \wt_i(b_2)=\vp_i(b_1 \ot
b_2)>\vep_i(b_3) - a_{ii}$. Since $\vep_i(b_3) - a_{ii} -
\wt_i(b_2) < \vp_i(b_1) \le \vep_i(b_2) - a_{ii}$, we have
$\vep_i(b_3) - \wt_i(b_2) < \vep_i(b_2)$, which implies $\vep_i
(b_2 \ot b_3) = \vep_i(b_2)$. Now we have $\vep_i (b_2 \ot b_3) <
\vp_i(b_1) \le \vep_i (b_2 \ot b_3) - a_{ii}$. Therefore we get
the desired result.


If $\vp_i (b_1) \le \vep_i (b_2)$, it suffices to show that
$\eit(b_1 \ot (b_2 \ot b_3)) = b_1 \ot (\eit b_2 \ot b_3)$. In
this case, we have $\vp_i(b_1 \ot b_2) = \vp_i (b_2)$. Hence by
our assumption, we obtain
\begin{align}\label{eq:1}
\vp_i(b_2) > \vep_i(b_3) - a_{ii}.
\end{align}
In particular, $\vep_i(b_2 \ot b_3) = \vep_i(b_2) \ge \vp_i(b_1)$,
which yields $\eit(b_1 \ot (b_2 \ot b_3)) = b_1 \ot \eit(b_2 \ot
b_3)$. By \eqref{eq:1}, we get what we wanted.

{\bf Case 2:} $\vep_i(b_3) < \vp_i (b_1 \ot b_2) \le \vep_i(b_3) -
a_{ii}$.

We will show $\eit (b_1 \ot (b_2 \ot b_3)) = 0$. By our
assumption, we have
\begin{align}\label{eq:2}
\vp_i(b_1) + \wt_i(b_2) \le \vep_i(b_3) - a_{ii} \quad \mbox{and}
\quad \vp_i(b_2) \le \vep_i(b_3) - a_{ii}.
\end{align}


If $\vp_i(b_2) \le \vep_i (b_3)$, we will show $\vp_i(b_1) +
a_{ii}\le\vep_i(b_2 \ot b_3)<\vp_i (b_1)$. Since $\vep_i(b_2 \ot
b_3) = \vep_i(b_3) - \wt_i(b_2)$, we must show $\vp_i(b_1) +
a_{ii}\le\vep_i(b_3) - \wt_i(b_2)<\vp_i (b_1)$. By \eqref{eq:2},
it suffices to show the second inequality. Since
$\vp_i(b_2)\le\vep_i(b_3)<\vp_i(b_1\otimes b_2)$, we have $\vp_i
(b_1 \ot b_2) = \vp_i (b_1) + \wt_i(b_2)$. Hence we obtain $\vp_i
(b_1) > \vep_i(b_3) - \wt_i(b_2)$.


If $\vep_i (b_3) < \vp_i(b_2)$, by the first inequality of
\eqref{eq:2}, we know $\vp_i (b_1) \le \vep_i(b_2 \ot b_3) -
a_{ii}$. Hence $\eit\bl(b_1 \ot (b_2 \ot b_3)\br) = 0$ or
$\eit\bl(b_1 \ot (b_2 \ot b_3)\br) = b_1 \ot \eit(b_2 \ot b_3)$.
By our assumption and the second inequality of \eqref{eq:2}, the
latter is also $0$.


{\bf Case 3:} $\vp_i (b_1 \ot b_2) \le \vep_i(b_3)$.

In this case, we have $\eit\bl((b_1\otimes b_2)\ot b_3\br)=
(b_1\otimes b_2)\ot \eit b_3$. Hence it is enough to show
$\vp_i(b_1) \le \vep_i(b_2 \ot b_3)$ and $\vp_i(b_2) \le
\vep_i(b_3)$. By the definition of $\vp_i$, we have (a)
$\vp_i(b_1) + \wt_i(b_2) \le \vep_i(b_3)$ and (b) $\vp_i(b_2) \le
\vep_i(b_3)$. Hence we have $\vp_i(b_1)\le
\vep_i(b_3)-\wt_i(b_2)\le\vep_i(b_2\otimes b_3)$, in which the
first inequality follows from (a). 
\end{proof}

\begin{remark}
The category of crystals $B$ such that
$$\text{$\vep_i(b)=0$, $\vp_i(b)=\wt_i(b)\ge0$
for any $i\in I^\im$ and $b\in B$}$$ is closed under tensor
product. Note that $B(\infty)$ and $B(\lambda)$ ($\la\in P^+$)
belong to this category.
\end{remark}

\begin{example} \label{exam:Tla}
For $\la \in P$, let $T_{\la}=\{t_{\la} \}$ and define
\begin{equation*}
\begin{aligned}
& \wt(t_{\la})=\la, \quad \eit t_{\la} = \fit t_{\la} =0 \quad
\text{for all} \ \ i\in I, \\
& \vep_i(t_{\la}) = \vphi_i(t_{\la}) = -\infty \quad \text{for
all} \ \ i\in I.
\end{aligned}
\end{equation*}
Then $T_{\la}$ is a crystal. We have $T_\la\otimes T_\mu\simeq
T_{\la+\mu}$. Note that $T_0$ is a unit object of the tensor
category of crystals (see e.g.\ \cite{KS}). Using this crystal,
Proposition~\ref{prop:projection} can be translated into following
the statement: \eq &&\text{for every $\la \in P^{+}$, there exists
an embedding $\pi_{\la}\cl B(\la)\to B(\infty)\otimes T_\la$.}
\eneq
\end{example}

\begin{example}
For each $i\in I$, let $B_i = \set{b_i(-n) }{ n\ge 0}$. Then $B_i$
is a crystal with the maps defined by
\begin{equation*}
\begin{aligned}
& \wt b_i(-n) = -n \alpha_i, \\
& \eit b_i(-n) = b_i(-n+1), \quad \fit b_i(-n) = b_i(-n-1), \\
& \tilde e_j b_i(-n) = \tilde f_j b_i(-n) = 0 \quad \text{if} \ \
j \neq i,\\
& \vep_i (b_i(-n)) = n, \quad \vphi_i (b_i(-n))=-n \quad \text{if}
\ \ i\in I^\re, \\
& \vep_i (b_i(-n)) = 0, \quad \vphi_i
(b_i(-n))=\wt_i(b_i(-n))=-na_{ii} \quad \text{if}
\ \ i\in I^\im, \\
& \vep_j (b_i(-n)) = \vphi_j (b_i(-n)) = -\infty \quad \text{if} \
j \neq i.
\end{aligned}
\end{equation*}
Here, we understand $b_i(-n)=0$ for $n<0$. The crystal $B_i$ is
called an {\it elementary crystal}.

\end{example}

\begin{example}\label{exam:tens}
For $\la,\mu\in P^{+}$, the tensor product $B(\la)\otimes B(\mu)$
is a crystal associated with $V(\la)\otimes V(\mu)$. There exists
a unique strict embedding:
$$\Phi_{\la,\mu}\cl B(\la+\mu)\To B(\la)\otimes B(\mu),$$
which sends $u_{\la+\mu}$ to $u_\la\otimes u_\mu$.
\end{example}

\begin{example}\label{exam:c}
Let $C=\{c\}$ be the crystal with $\wt(c)=0$ and
$\vep_i(c)=\vp_i(c)=0$, $\fit c=\eit c=0$ for any $i\in I$. Then
$C$ is isomorphic to $B(0)$. For a crystal $B$, $b\in B$ and $i\in
I$, we have \eqn
\wt(b\otimes c)&=&\wt(b),\\
\vep_i(b\otimes c)&=&\max(\vep_i(b),-\wt_i b),\\
\vp_i(b\otimes c)&=&\max(\vp_i(b),0),\\
\eit(b\otimes c)&=&
\begin{cases}
\eit b\otimes c&\text{if $\vp_i(b)\ge0$ and $i\in I^\re$,}\\
\eit b\otimes c&\text{if $\vp_i(b)+a_{ii}>0$ and $i\in I^\im$,}\\
0&\text{otherwise,}
\end{cases}\\
\fit(b\otimes c)&=&
\begin{cases}
\fit b\otimes c&\text{if $\vp_i(b)>0$,}\\
0&\text{otherwise.}
\end{cases}
\eneqn

In general, $B \ot C$ is {\emph not} isomorphic to $B$.
\end{example}

\begin{example} \label{exam:B(iota)}
Let $\bfi=(i_1,i_2,\dots)$ be an infinite sequence in $I$ such
that every $i\in I$ appears infinitely many times in $\bfi$. For
$k\in\Z_{>0}$, set $B(k)=B_{i_{k}}\otimes\cdots \otimes B_{i_1}$.
For $k_1\le k_2$, let $\psi_{k_2,k_1}\cl B(k_1) \to B(k_2)$ be the
map $b\mapsto b_{i_{k_2}}(0)\otimes\cdots\otimes
b_{i_{k_1+1}}(0)\otimes b$. Then  $\{B(k)\}_{k\ge 1}$ is an
inductive system. We shall consider the (set-theoretical)
inductive limit of $B(k)$: \eqn&& B({\bfi})=\{\cdots\otimes
b_{i_k}(-x_k)\otimes \cdots\otimes b_{i_1}(-x_1)\\
&&\hs{13ex}\in \cdots\otimes B_{i_k}\otimes \cdots \otimes B_{i_1}
\,;\, \text{$x_k\in \Z_{\ge 0}$, and $x_k=0$ for $k\gg 0$}\}.
\eneqn Let $\psi_k\cl B(k)\to B({\bfi})$ be the canonical
injective map. The maps $\psi_{k_2,k_1}$ are not crystal
morphisms, but they have the following properties: \bnum \item
$\wt$ is preserved by $\psi_{k_2,k_1}$, \item for $i\in I$, $k\in
\Z_{>0}$ and $b\in B(k)$, the sequences
$\{\psi_{k'}(\eit(\psi_{k',k}(b))\}_{k'\ge k}$,
$\{\psi_{k'}(\fit(\psi_{k',k}(b))\}_{k'\ge k}$ and
$\{\vep_i(\psi_{k',k}(b))\}_{k'\ge k}$,
$\{\vphi_i(\psi_{k',k}(b))\}_{k'\ge k}$ are stationary, \ee
Therefore, the inductive limit $B({\bfi})$ has a crystal
structure. This is explicitly given as follows. Let $b=\cdots
\otimes b_{i_k}(-x_k)\otimes \cdots\otimes b_{i_1}(-x_1)\in
B({\bfi})$. Then we have \eqn \wt(b)&=&-\sum_k x_k\alpha_{i_k}.
\eneqn For $i\in I^\re$, we have
\begin{equation*}\aligned
\vep_i(b)&=\max\Bigl\{x_{k}+\sum_{l>k}\langle
h_i,\alpha_{i_l}\rangle x_{l}\,;\,1\le k,\,i=i_k\Bigr\},\\
\vp_i(b)&=\max \Bigl\{-x_k-\sum_{1\le l<k}\langle
h_i,\alpha_{i_l}\rangle x_{l}\,;\,1\le k,\, i=i_k\Bigr\},
\endaligned
\end{equation*}
and, for $i\in I^\im$, we have
$$\vep_i(b)=0\quad\text{and}\quad \vp_i(b)=\wt_i(b).$$
For $i\in I^\re$, we have
\begin{equation*}\aligned
\eit b&=\begin{cases} \cdots \otimes
b_{i_{n_e}}(-x_{n_e}+1)\otimes\cdots
\otimes b_{i_1}(-x_1)&\text{if $\vep_i(b)> 0$,}\\[1ex]
0&\text{if $\vep_i(b)\le 0$,}
\end{cases}\\
\fit b&= \cdots \otimes b_{i_{n_f}}(-x_{n_f}-1)\otimes\cdots
\otimes b_{i_1}(-x_1),
\endaligned
\end{equation*}
where $n_e$ (resp.\ $n_f$) is the largest (resp.\ smallest)
$k\ge1$ such that $i_k=i$ and $x_{k}+\sum_{l>k}\langle
h_i,\alpha_{i_l}\rangle x_{l}=\vep_i(b)$. Note that such an $n_e$
exists if $\vep_i(b)>0$.

When $i\in I^\im$, let $n_f$ be  the smallest $k$ such that
\begin{equation*}
i_k=i \quad\text{and}\quad \sum_{l>k}\langle
h_i,\alpha_{i_l}\rangle x_l=0.
\end{equation*}
Then we have
\begin{equation*}
\fit b=\dots\otimes b_{i_{n_f}}(-x_{n_f}-1)\otimes\dots\otimes
b_{i_1}(-x_1)
\end{equation*}
and
\begin{equation*}
\eit b=\left\{ \ba{l} \dots\otimes
b_{i_{n_f}}(-x_{n_f}+1)\otimes\dots\otimes b_{i_1}(-x_1)\\[1ex]
\hs{6ex}\parbox{45ex}{if
$x_{n_f}>0$ and $\sum_{k<l\le n_f}\lan h_i,\alpha_{i_l}\ran
x_l<a_{ii}$ for any $k$
such that $1\le k<n_f$ and $i_{k}=i$,}\\[3ex]
0\qquad\text{otherwise.} \ea\right.
\end{equation*}

\end{example}

\section{Crystal Embedding Theorem}

In this section, we will prove one of the main results of this
paper, the {\it crystal embedding theorem} for  quantum
generalized Kac-Moody algebras.

\begin{thm}
For all $i\in I$, there exists a unique strict embedding
$$\Psi_i\cl B(\infty)\to B(\infty)\otimes B_i,$$
called the crystal embedding.
\end{thm}
The map $\Psi_i$ sends ${\mathbf{ 1}}$ to ${\mathbf{1}}\otimes
b_i(0)$, because there is a unique vector of weight $0$ in
$B(\infty)\otimes B_i$.

\begin{proof}
Let $b=\tilde{f}_{i_1}\cdots\tilde{f}_{i_r} \mathbf{1}\in
B(\infty)$. Take $\mu\gg 0$ in $P^+$ such that $b\in
\Im(\pi_{\mu})$, where $\pi_{\mu}\cl B(\mu)\rightarrow B(\infty)$
is the map given in Proposition \ref{prop:projection}. Hence
$b_\mu\seteq \tilde{f}_{i_1}\cdots\tilde{f}_{i_r} u_\mu\in B(\mu)$
satisfies $\pi_\mu(b_\mu)=b$.

Set $l=\mu(h_i)$ and set $\la=\mu-l\Lambda_i\in P^+$.
Then $\la(h_i)=0$ and there is a unique strict embedding
$\Phi_{\la,l\La_i}\cl B(\mu)\rightarrow B(\la)\otimes B(l\La_i)$,
which sends $u_{\mu}$ to $u_{\la}\otimes u_{l\La_i}$ (see
Example~\ref{exam:tens}). We claim that \eq &&\quad\parbox{65ex}{
$\Phi_{\la,l\La_i}(b_{\mu})=
b'\otimes \fit^n u_{l\La_i}$ for some $b'$ and $n\in\Z_{\ge0}$,\\
moreover $\pi_\lambda(b')\otimes b_i(-n)$ does not depend on the
choice of $\mu\gg0$.} \label{eq:claim} \eneq Then we define
$\Psi_i(b)=\pi_\lambda(b')\otimes b_i(-n)$.

We show \eqref{eq:claim} by the induction on $r$. If $r=0$, our
assertion is
obvious. 
Assume that our assertion is true for $r-1$ and let
$b_1=\tilde{f}_{i_2}\cdots\tilde{f}_{i_r}u_\mu$. By the induction
hypothesis, we have
\begin{equation*}
\aligned \Phi_{\la,l\La_i} (b_1)&=b'_1\otimes \fit^m u_{l\La_i}
\endaligned
\end{equation*}
for some $b'_1\in B(\lambda)$ and $m\in\Z_{\ge0}$. Hence it
suffices to show that
\begin{equation}\label{eq:b1i}
\parbox{60ex}{$\tilde{f}_{i_1}(b'_1\otimes \fit^m u_{l\La_i})
=b'\otimes \fit^nu_{l\La_i}$ for some $b'\in B(\la)$
and $n\in\Z_{\ge0}$,\\
and $\tf_{i_1}(\pi_\la(b'_1)\otimes b_i(-m))= \pi_\la(b')\otimes
b_i(-n)$.}
\end{equation}
If $i_1=i$, then $\vp_i(b'_1)=\vp_i(\pi_\la(b'_1))$ since
$\la(h_i)=0$, and
$$\vep_i(\fit^mu_{l\Lambda_i})=\vep_i(b_i(-m))
=\begin{cases}
m&\text{if $i\in I^\re$,}\\
0&\text{if $i\in I^\im$,}\\
\end{cases}
$$
and hence \eqref{eq:b1i} is obvious. Suppose that $i_1\neq i$.
Then
\begin{equation*}
\vp_{i_1}(b'_1)=\vp_{i_1}(\pi_{\la}(b'_1)) +\lan
h_{i_1},\la\rangle \gg 0=\vep_{i_1}(\fit^m u_{l\La_i}),
\end{equation*}
we have
\begin{equation*}
\tf_{i_1}(b'_1\otimes \fit^m u_{l\La_i})=\tf_{i_1}b'_1\otimes
\fit^m u_{l\La_i}.
\end{equation*}
On the other hand we have $\tf_{i_1}(\pi_\la(b'_1)\otimes
b_i(-m))= \tf_{i_1}\pi_\la(b'_1)\otimes b_i(-m)$, which proves our
claim \eqref{eq:b1i}.


It is straightforward to verify that
%
$\Psi_i\cl B(\infty)\rightarrow B(\infty)\otimes B_i$ is a strict
crystal embedding. The uniqueness of $\Psi_i$ is obvious since
$B(\infty)\otimes B_i$ has a unique vector with weight $0$.
\end{proof}

The crystal embedding theorem yields a procedure to determine the
structure of the crystal $B(\infty)$ in terms of elementary
crystals. Take an infinite sequence ${\bfi} =(i_1,i_2,\dots)$ in
$I$ such that every $i\in I$ appears infinitely many times. Such a
sequence always exists since $I$ is countable. For each $N\ge 1$,
taking the composition of crystal embeddings repeatedly, we obtain
a strict crystal embedding
\begin{equation}
\aligned &\Psi^{(N)}\seteq(\Psi_{i_N}\otimes \id\otimes
\cdots\otimes
\id)\circ\cdots\circ(\Psi_{i_2}\otimes \id)\circ \Psi_{i_1}\cl \\
&B(\infty)\hookrightarrow B(\infty)\otimes B_{i_1}\hookrightarrow
B(\infty)\otimes B_{i_2}\otimes B_{i_1}\hookrightarrow\\
&\qquad\qquad\qquad\qquad\qquad\qquad \cdots \hookrightarrow
B(\infty)\otimes B_{i_N}\otimes \cdots\otimes B_{i_1}.
\endaligned
\end{equation}
It is easily seen that, for any $b\in B$, there exists $N>0$ such
that
\begin{equation*}
\Psi^{(N)}(b)=\mathbf{1}\otimes b_{i_N}(-x_N)\otimes \cdots\otimes
b_{i_1}(-x_1)
\end{equation*}
for some $x_1,\dots,x_N\in \Z_{\ge 0}$. Set $x_k=0$ for $k>N$.
Then for any $k\ge N$, we have $\Psi^{(k)}(b)=\mathbf{1}\otimes
b_{i_k}(-x_k)\otimes \cdots\otimes b_{i_1}(-x_1)$. Hence,
associating $\cdots\ot b_{i_{N+1}}(0) \otimes b_{i_N}(-x_N)\otimes
\cdots\otimes b_{i_1}(-x_1)$ to $b$, we obtain a map $B(\infty)\to
B(\mathbf{i})$ (see Example~\ref{exam:B(iota)}). We can easily see
that it is a crystal morphism, and we obtain the following result.


\begin{prop}
$B(\infty)$ is strictly embedded in the crystal $B({\bfi})$
introduced in {\rm Example~\ref{exam:B(iota)}}.
\end{prop}
Hence $B(\infty)$ is isomorphic to the connected component of
$B({\bfi})$ containing $b(\bfi,\mathbf{0})\seteq\cdots\otimes
b_{i_2}(0)\otimes  b_{i_1}(0)$.

\begin{example}\label{exam:rank2}

Let $I=\{1,2\}$ and consider the quantum generalized Kac-Moody
algebra $U_q(\g)$ associated with a rank 2 Borcherds-Cartan matrix
$$A= \left( \begin{matrix} 2 & -a \\ -b & -c \end{matrix} \right)
\quad \text{for some $a, b \in \Z_{>0}$ and $c \in 2\Z_{\ge
0}$}.$$ Take an infinite sequence ${\bfi}=(1,2,1,2,\dots)$ and let
$B$ be the set of elements of the form
$$b({\bfi}, \mathbf{x}) \seteq\cdots\otimes
b_{2}(-x_{2k})\otimes b_1(-x_{2k-1})\otimes \cdots\otimes
b_2(-x_2)\otimes b_{1}(-x_1)\in B({\bfi})$$ satisfying the
following conditions: \bnum
\item $ax_{2k}-x_{2k+1}\ge 0$ for all $k\ge 1$,

\item for all $k\ge 2$ such that $x_{2k}>0$, we have $x_{2k-1}>0$
and $ax_{2k}-x_{2k+1}> 0$. \ee It is shown in \cite{Shin05} that
$B$ is the connected component of $B({\bfi})$ containing
$b({\bfi}, \mathbf{0})\seteq \cdots\otimes b_{2}(0)\otimes
b_{1}(0)\otimes b_2(0)\otimes b_{1}(0)$. Therefore, $B$ is
isomorphic to the crystal $B(\infty)$.
\end{example}

\begin{example} \label{exam:Monster} 

Let $I=\set{(i,t)}{i\in\Z_{\ge -1},\,1\le t\le c(i)}$, where
$c(i)$ is the $i$-th coefficient of the elliptic modular function
$$j(q)- 744 = q^{-1} + 196884q + 21493760q^2 + \cdots
= \sum_{i=-1}^{\infty} c(i) q^i. $$ Consider the Borcherds-Cartan
matrix $A=(a_{(i,t),(j,s)})_{(i,t),(j,s) \in I}$  whose entries
are given by $a_{(i,t),(j,s)}=-(i+j)$.
The associated generalized Kac-Moody algebra $\mathfrak g$ is
called the {\it Monster Lie algebra}, and it played a crucial role
in Borcherds' proof of the Moonshine conjecture \cite{Bo92}. More
precisely, Borcherds derived the {\it twisted denominator
identity} for the Monster Lie algebra with the action of the
Monster, from which the replication formulae for the Thompson
series follow.

In this paper, we deal with the corresponding quantum group
$U_q(\mathfrak g)$ which we call the {\it quantum Monster
algebra}. Take the infinite sequence
\begin{equation*}
\begin{aligned}
{\bf i} = ({\bf i}(k))_{k=1}^{\infty}=(&
(-1,1),(1,1),\dots,(1,c(1)); (-1,1), (1,1),\dots,
(1,c(1)),\\
&(2,1),\dots,(2,c(2)); (-1,1), (1,1),\dots, (1,c(1)),(2,1),
\dots, \\
& (2,c(2)),(3,1),\dots,(3,c(3));(-1, 1), \dots).
\end{aligned}
\end{equation*}
Note that $(-1,1)$ appears at the $b(n)$-th position for $n\ge 0$,
where
\begin{equation*}
b(n)=nc(1)+(n-1)c(2)+\dots+c(n)+n+1.
\end{equation*}
For $k\in\Z_{>0}$, we denote by $k^{(-)}$ the largest integer
$l<k$ such that ${\bf i}(l)={\bf i}(k)$. If such an $l$ does not
exist, then set $k^{(-)}=0$. Let $B$ be the set of elements
\begin{equation*}
b({\bf i}, \mathbf{x})= \cdots\otimes b_{{\bf
i}(k)}(-x_{k})\otimes\cdots\otimes b_{{\bf i}(1)}(-x_1)\in B({\bf
i})
\end{equation*}
satisfying the following conditions: \bnum
\item $x_{b(1)}=0$,

\item for all $n\ge 1$, we have
\begin{equation*}
\aligned
-\hs{-3ex}\sum_{b(n)<l<b(n+1)}\lan
h_{(-1,1)},\alpha_{\mathbf{i}(l)}\ran x_l \ge x_{b(n+1)},
\endaligned
\end{equation*}

\item if ${\bf i}(k)\neq (-1,1)$, $x_k>0$ and
$k^{(-)}>0$, then
$$\sum_{k^{(-)}<l<k}\langle h_{{\bf i}(k)}, \alpha_{{\bf i}(l)}\rangle x_l<
0.$$ In addition, if $x_l=0$ for all $k^{(-)}<l<k$ such that ${\bf
i}(l)\neq (-1,1)$, then we have
\begin{equation*}
\aligned
-\hs{-3ex}\sum_{b(n)<l<b(n+1)}\lan
h_{(-1,1)},\alpha_{\mathbf{i}(l)}\ran x_l
>x_{b(n+1)},
\endaligned
\end{equation*}
where $n$ is a unique integer such that $k^{(-)}<b(n)<k$. \ee Then
$B$ is the conncected component of $B({\bfi})$ containing
$b(\mathbf{i}, \mathbf{0})= \cdots\otimes
b_{\mathbf{i}(2)}(0)\otimes b_{{\bf i}(1)}(0)$ (see
\cite{Shin05}). Therefore, $B$ is isomorphic to the crystal
$B(\infty)$.

\end{example}

\vskip 1cm
\section{Characterization of $B(\infty)$ and $B(\la)$}

As an application of the crystal embedding theorem, we will give a
characterization of the crystals $B(\infty)$ and $B(\la)$.

\begin{thm} \label{thm:B(infty)}
Let $B$ be a crystal. Suppose that $B$ satisfies the following
conditions: \bnum
\item  $\wt(B)\subset -Q_{+}$,

\item there exists an element $b_0\in B$ such that
$\wt(b_0)=0$,


\item for any $b\in B$ such that $b\neq b_0$, there exists
some $i\in I$ such that $\eit b\neq 0$,

\item for all $i$, there exists a strict embedding
$\Psi_i\cl B\rightarrow B\otimes B_i$. \ee Then there is a crystal
isomorphism
$$B\isoto B(\infty),$$
which sends $b_0$ to $\mathbf{1}$.
\end{thm}

\begin{proof}
Note that $b_0$ is a unique element of weight $0$ in $B$. Indeed,
if $b\not=b_0$ has weight $0$, then taking $i\in I$ such that
$\eit b\in B$, the weight of  $\eit b$ does not belong to $-Q_+$,
which is a contradiction.

Note also that $\eit b_0=0$ for any $i\in I$. Indeed, otherwise,
$\wt(\eit b_0)\not\in-Q_+$.

\smallskip
Since $B\otimes B_i$ has a unique vector of weight $0$, we have
$\Psi_i(b_0)=b_0\otimes b_i(0)$.

Take an infinite sequence ${\bfi}=(i_1,i_2,\dots)$ in $I$ such
that every $i\in I$ appears infinitely many times.
Similarly to the case of $B(\infty)$, we obtain the strict crystal
embeddings
\begin{equation}
\ba{l} \Psi^{(N)}\seteq(\Psi_{i_N}\otimes \id\otimes \cdots\otimes
\id)\circ\cdots\circ(\Psi_{i_2}\otimes \id)\circ \Psi_{i_1}\cl \\
B(\infty)\hookrightarrow B(\infty)\otimes B_{i_1}\hookrightarrow
B(\infty)\otimes B_{i_2}\otimes B_{i_1}\hookrightarrow \\
\qquad\qquad\qquad\qquad\qquad\qquad\cdots \hookrightarrow
B(\infty)\otimes B_{i_N}\otimes \cdots\otimes B_{i_1}. \ea
\end{equation}
We can easily see that for any $b\in B$, there exists $N\gg0$ such
that
\begin{equation*}
\aligned \Psi^{(N)}(b)
%
&=b_0\otimes b_{i_N}(-x_N)\otimes \cdots\otimes b_{i_1}(-x_1)
\endaligned
\end{equation*}
for some $x_1,\dots,x_N\in \Z_{\ge 0}$. Hence we get a strict
crystal embedding $B\hookrightarrow B({\bfi})$ given by
\begin{center}
$b\mapsto 
\cdots\otimes b_{i_{N+1}}(0)\otimes b_{i_{N}}(-x_N)\otimes
\cdots\otimes b_{i_1}(-x_1)$.
\end{center}
Since this embedding sends $b_0$ to
$$b({\bfi}, \mathbf{0}) \seteq \cdots\otimes
b_{i_{N+1}}(0)\otimes b_{i_N}(0)\otimes \cdots\otimes
b_{i_1}(0),$$ $B$ is isomorphic to the the connected component of
$B({\bfi})$ containing $b({\bfi}, {\mathbf{0}})$.
The theorem then follows from the fact that $B(\infty)$ is also
isomorphic to the smallest full subcrystal of $B({\bfi})$
containing $b({\bfi}, \mathbf{0})$.
\end{proof}

For $\la\in P^+$, the crystal $B(\la)$ is embedded in
$B(\infty)\otimes T_\lambda$.

\begin{thm} \label{thm:B(la)}
Let $\la\in P^+$ be a dominant integral weight. Then $B(\lambda)$
is isomorphic to the connected component of $B(\infty)\otimes
T_\la\otimes C$ containing $\mathbf{1}\otimes t_\lambda\otimes c$.
Here $C$ is the crystal introduced in {\rm Example \ref{exam:c}}.
%
%
%
%
%
%
%
%
\end{thm}
\begin{proof}
Note that the embedding $\iota\cl B(\la)\to B(\infty)\otimes
T_\la$ commutes with $\eit$. Let us remark that for $b\in B(\la)$,
$\fit(b)=0$ if and only if $\vp_i(b)\le 0$. Hence, $\iota$ induces
a strict embedding $B(\la)\to B(\infty)\otimes T_\la\otimes C$.
Hence the assertion follows.
\end{proof}

\begin{example}
Let $U_q(\g)$ be the quantum generalized Kac-Moody algebra
associated with the rank 2 Borcherds-Cartan matrix given in
Example \ref{exam:rank2}. We shall use the notations there. Let
$\la \in P^{+}$. Then $B(\la)$ is isomorphic to the connected
component $B^\la$ of $B({\bfi})\otimes T_{\la}\otimes C$
containing $b({\bfi}, \mathbf{0}) \ot t_{\la}\ot c$. It was shown
in \cite{Shin06} that $B^\la$ is the set of elements of the form
$$ b({\bfi}, \mathbf{x}) \ot t_{\la}\ot c
\seteq\cdots\otimes b_{2}(-x_{2k})\otimes b_1(-x_{2k-1})\otimes
\cdots\otimes b_2(-x_2)\otimes b_{1}(-x_1)\otimes t_{\la}\ot c$$
satisfying the conditions (i) and (ii) in Example~\ref{exam:rank2}
and two additional conditions: \be[\quad (a)]

\item $0\le x_1\le \langle h_1,\la\rangle$,

\item if $x_2> 0$ and $\lan h_2, \la \ran =0$, then $x_1 >0$.
\ee


\end{example}

\begin{example}
Let $U_q(\g)$ be the quantum Monster algebra in Example
\ref{exam:Monster}, and let $B(\la)$ be the irreducible highest
weight crystal with $\la \in P^{+}$. Using the notations in
Example \ref{exam:Monster}, it was shown in \cite{Shin06} that the
connected component of $B({\bfi})\otimes T_{\la}\ot C$ containing
$b({\bfi}, \mathbf{0}) \ot t_{\la}\ot c = \cdots\otimes b_{{\bf
i}(k)}(0)\otimes \cdots\otimes b_{{\bf i}(1)}(0)\otimes t_{\la}\ot
c$ is the set $B^{\la}$ consisting of elements of the form
\begin{equation*}
b({\bf i}, \mathbf{x}) \ot t_{\la}\ot c = \cdots\otimes b_{{\bf
i}(k)}(-x_{k})\otimes \cdots \otimes b_{{\bf i}(1)}(-x_{1})\otimes
t_{\la}\ot c
\end{equation*}
satisfying the conditions (i)--(iii) in Example~\ref{exam:Monster}
and two additional conditions: \be[\quad (a)]
\item $0\le x_1\le \langle h_{(-1,1)},\la\rangle$,

\item if $\mathbf{i}(k)\neq (-1,1)$, $\lan h_{\mathbf{i}(k)},\la\ran=0$,
$x_k>0$ and $k^{(-)}=0$, then there exists $l$  such that $1\le
l<k$, $\lan h_{\mathbf{i}(k)},\alpha_{\mathbf{i}(l)}\ran<0$ and
$x_l>0$.
\ee

Hence by Theorem \ref{thm:B(la)}, we conclude that $B^{\la}$ is
isomorphic to the crystal $B(\la)$.
\end{example}


\providecommand{\bysame}{\leavevmode\hbox
to3em{\hrulefill}\thinspace}

\end{document}

--